\documentclass[a4paper,11pt]{article}
\pdfoutput=1 

\usepackage{mathrsfs}
\pdfoutput=1 

\usepackage{jheppub} 

\usepackage[T1]{fontenc}

\usepackage{mathrsfs}
\usepackage{stmaryrd}
\usepackage{epsfig}
\usepackage{graphicx}
\usepackage{graphics}
\usepackage{latexsym}
\usepackage{amssymb,amsthm,amsfonts,amsmath,amscd}
\usepackage[arrow, matrix, curve]{xy}
\usepackage{syntonly}
\usepackage{extarrows}
\usepackage{txfonts}

\theoremstyle{plain}
\newtheorem{thm}{Theorem}[section]
\newtheorem{theorem}[thm]{Theorem}
\newtheorem{lemma}[thm]{Lemma}

\newtheorem{proposition}[thm]{Proposition}

\theoremstyle{definition}

\newtheorem{definition}[thm]{Definition}

\newtheorem{example}[thm]{Example}

\title{HOMFLY polynomial from a generalized Yang-Yang function}

\author{Hu Sen}
\author{and Liu Peng}

\affiliation{School of Mathematics, University of Science and Technology of China}
\affiliation{Wu Wen-Tsun Key Lab of Mathematics, Chinese Academy of Sciences}

\emailAdd{shu@ustc.edu.cn}
\emailAdd{pliu@mail.ustc.edu.cn}

\abstract{Starting from the free field realization of Kac-Moody Lie algebra, we define a generalized Yang-Yang function. Then for the Lie algebra of type $A_{n}$, we derive braiding and fusion matrix by braiding the thimble from the generalized Yang-Yang function. One can construct a knots invariant $H(K)$ from the braiding and fusion matrix. It is an isotropy invariant and obeys a skein relation. From them, we show that the corresponding knots invariant is HOMFLY polynomial.}

\begin{document}
\maketitle
\flushbottom

\section{Introduction}
\label{sec:intro}

Knots invariants are topological invariants under Reidemeister moves, among which Jones polynomial~\cite{a} and its generalization HOMFLY polynomial~\cite{b} have been well studied. Statistical mechanics explanation~\cite{c} of knots invariants was given and also the quantum field theory method has been used to recover the Jones polynomial as a partition function (or an expectation value) for given knots~\cite{d}. Recently, Gaiotto and Witten~\cite{e} reconstruct Jones polynomial from the complexified Chern-Simons theory by studying opers structure, integrable Gaudin model and Virasoro conformal block in conformal field theory. The conformal block that comes from the integral of the Chern-Simons functional over an infinite dimensional thimble is shown to be equivalent to that from the integral of a Yang-Yang function over a finite dimensional thimble. It gives a powerful tool to reveal the relationship between the Chern-Simons gauge theory and knots invariants. In this paper, we generalize the method in Gaiotto and Witten's work~\cite{e} to study the thimble of the generalized Yang-Yang function associated to a simple Lie algebra and derive braiding matrix and skein relation for the Lie algebra of type $A_{n}$.

In section 2, we give a brief review of the method in~\cite{e}. In section 3, we use Wakimoto's free field realization of Kac-Moody algebra to construct a generalized Yang-Yang function as our start point, then introduce the concept of thimble and derive the braiding and fusion matrix from braiding thimbles of the generalized Yang-Yang function associated with $A_{n}$ Lie algebra. Finally, framing independent knots invariant is constructed and its skein relation is derived, which shows that it is HOMFLY polynomial.

Based on works on Landau-Ginzburg B models ~\cite{j,k}, we ~\cite{l} derive braiding from the $t t^{*}$ equations of Landau-Ginzburg B models with Yang-Yang function the super-potential function.

Knot invariants for the $B$, $C$, $D$ type Lie algebra correspond to Kauffman polynomials, details will appear elsewhere ~\cite{m}.

\section{A brief review of Gaiotto-Witten method}

In~\cite{e}, the relationship between the complexified Chern-Simons gauge theory and Jones polynomial was studied. The gradient flow of the complexified Chern-Simons functional has a nice structure: with the boundary condition incorporating the information of knots at the finite boundary and symmetry breaking at the infinity, solutions of the gradient flow equation are corresponding to opers with monodromy free singularities.

When the gauge group is $SL(2,\mathbb{C})$, an oper with monodromy free singularities satisfies the equation:
\begin{equation}
\partial z\dfrac{P(z)}{Q(z)}=-\dfrac{K(z)}{Q^{2}(z)},
\end{equation}
where $K(z)=\prod^{d}_{a=1}(z-z_{a})^{\lambda_{a}}$ encodes the position and the charge of the oper singularities and $Q(z)=\prod^{\mathcal {q}}_{j=1}(z-w_{j})$ is one component of the section $\left(
                                                                              \begin{array}{c}
                                                                                P(z) \\
                                                                                Q(z) \\
                                                                              \end{array}
                                                                            \right)
$.
The residue free condition of the left hand side leads to
\begin{equation}
\label{eq:a}
\sum _{a}\dfrac {\lambda_{a}} {w_{j}-z_{a}}=\sum _{s\neq j}\dfrac {2} {w_{j}-w_{s}}, \quad j=1,2...,p.
\end{equation}
This equation is called the Bethe equation and its solutions $w_{j}$ Bethe roots. It can be considered as critical point equation \begin{equation}
\frac{\partial \mathscr{W}(w,z)}{\partial w_{j}}=0,\quad j=1,2,...,p
\end{equation} of the Yang-Yang function
\begin{equation}
\label{eq:b}
\mathscr{W}(w,z)=\sum _{j,a}\lambda _{a}\ln(w_{j}-z_{a}) -\sum _{s<j}2\ln(w_{j}-w_{s})-\sum _{a<b}\frac{\lambda_{a}\lambda_{b}}{2}\ln(z_{a}-z_{b}).
\end{equation}
In~\cite{e}, the free field realization of Virasoro conformal blocks was used to construct the representation of braid group. It gives an integral formula
\begin{multline}
\int _{\Gamma }\langle \prod _{i}V_{1/b}\left( w_{i}\right) \prod _{a}V_{-\lambda_{a} / 2b}\left( z_{a}\right) \rangle_{\mathrm{free}} \prod _{i}dw_{i}
\\=\int _{\Gamma }\prod_{i,a}(w_{i}-z_{a})^{\frac{\lambda_{a}}{b^{2}}}\prod_{i<j}(w_{i}-w_{j})^{-\frac{2}{b^{2}}}\prod_{a<b}(z_{a}-z_{b})^{-\frac{\lambda_{a}\lambda_{b}}{2b^{2}}}\prod _{i}dw_{i}
=\int _{\Gamma }e^{\frac{\mathscr{W}}{b^{2}}}\prod _{i}dw_{i},
\end{multline}
 where $\Gamma$ is a thimble of the Yang-Yang function $\mathscr{W}$ and $b^{2}=-(k+2)$ is a constant relative to the level $k$ of the Chern-Simons theory. The definition of the thimble will be given in Section 3.
From conformal field theory, Virasoro conformal block, the solution of the Knizhnik-Zamolodchikov equation, gives a representation of braid group. So we can use braiding of the thimble of the Yang-Yang function to study the representation of braid group. As a multiple valued function, when thimble is braided, $e^{\frac{\mathscr{W}}{b^{2}}}$ integrated over the thimble will produce an additional phase factor. As is shown in the~\cite{e}, we can use thimbles as bases to represent braid group as a braiding matrix. Combining with the fusion matrix in a standard way~\cite{c} gives the famous Jones polynomial.

\section{Generalized Yang-Yang function and braiding of thimble}

\subsection{Generalized Yang-Yang function}

First, we make an introduction to Gaudin model~\cite{f} associated to a finite dimensional complex simple Lie algebra $g$ of rank $r$. $\Pi =\{\alpha_{1},\alpha_{2},...,\alpha_{r}\}$ is the set of the simple roots of $g$, $\{F_{i}, G_{i}, H_{i}\}, i=1,2,...,r$ the standard generators of $g$, $\{T_{\beta }\}$  the bases of $g$ and $\{T^{\beta}\}$ the dual bases induced by the Killing form of $g$. Considering distinct points $z_{1},...,z_{d}$ in $\mathbb{C}$, we associate each point an irreducible highest weight representation $V_{\lambda_{a}}$ of $g$, where $\lambda_{a}$ is a dominant integral weight. Thus $V_{\lambda_{a}}$ is a finite dimensional irreducible highest weight representation of $g$. $V_{(\lambda_{a})}\triangleq V_{\lambda_{1}}\otimes V_{\lambda_{2}}\otimes...\otimes V_{\lambda_{d}} $. The Hamiltonian operator is defined as
\begin{equation}
\Xi _{a}\triangleq \sum_{b\neq a}\dfrac{\sum_{\beta}T^{(a)}_{\beta}T^{(b)\beta}}{z_{a}-z_{b}}, \quad a=1,2,...,d
\end{equation}
 on $V_{(\lambda_{a})}$, where $T^{(a)}_{i}$ only acts on the component $V_{\lambda_{a}}$ of $V_{(\lambda_{a})}$. $v_{\lambda_{a}}$ is the highest weight vector of $V_{\lambda_{a}}$, then it is clear that $v_{\lambda_{1}}\otimes v_{\lambda_{2}}\otimes ... v_{\lambda_{d-1}}\otimes v_{\lambda_{d}}$ is a joint eigenvector of $\Xi _{a}$. To find other joint eigenvectors of $\Xi_{a}$ leads to the Bethe ansatz (for more detail see~\cite{f}). $w_{j} (j=1,2,...,p)$ are distinct points on $\mathbb{C}$ different from $z_{a}$. Each $w_{j}$ is associated with a colour $\alpha_{i_{j}}$, a simple root of $g$, where $i_{j}\in \{1,2,...,r\}$. Wakimoto realization is known as the free field realization of an affine Kac-Moody Lie algebra at the arbitrary level. In~\cite{f,g,i}, it was used to study the Bethe ansatz. The general Bethe equation for the simple Lie algebra $g$ is obtained:
\begin{equation}
\label{eq:c}
\sum _{a}\dfrac {(w_{i_{j}},\lambda_{a})} {w_{j}-z_{a}}=\sum _{s\neq j}\dfrac {(\alpha_{i_{j}},\alpha_{i_{s}})} {w_{j}-w_{s}}, \quad j=1,2,...,p,
\end{equation}
where $(,)$ is the inner product on the weight space induced by the Killing form of the simple Lie algebra $g$.
The correlation function of Wakimoto realization at arbitrary level $k$ gives conformal blocks for WZW model and a representation of braid group. It is a generalization of Virasoro conformal blocks:
\begin{equation}
\int _{\Gamma }\prod_{j,a}(w_{j}-z_{a})^{-\frac{(\alpha_{i_{j}},\lambda_{a})}{k+h^{\vee }}}\prod_{j<s}(w_{j}-w_{s})^{\frac{(\alpha_{i_{j}},\alpha_{i_{s}})}{k+h^{\vee }}}\prod_{a<b}(z_{a}-z_{b})^{\frac{(\lambda_{a},\lambda_{b})}{k+h^{\vee }}}\prod _{j}dw_{j}=\int _{\Gamma }e^{-\frac{\mathscr{W}}{k+h^{\vee }}}\prod _{j}dw_{j}
\end{equation}
with condition that
\begin{equation}
\sum_{\alpha}\lambda_{\alpha}-\lambda_{\infty}=\sum_{j}\alpha_{i_{j}},
\end{equation}
where $h^{\vee}$ is the dual Coxeter number.
We call
\begin{equation}
\mathscr{W}(w,z)=\sum _{j,a}(\alpha_{i_{j}},\lambda _{a})\ln(w_{j}-z_{a}) -\sum _{s< j}(\alpha_{i_{j}},\alpha_{i_{s}})\ln(w_{j}-w_{s})-\sum _{a< b}(\lambda_{a},\lambda_{b})\ln(z_{a}-z_{b})
\end{equation}
a generalized Yang-Yang function. It is obvious that the critical point equation of the generalized Yang-Yang function is the general Bethe equation \eqref{eq:c}.

With the symmetry breaking~\cite{e}, the generalized Yang-Yang function becomes:
\begin{multline}
\mathscr{W}(w,z)
=\sum _{j,a}(\alpha_{i_{j}},\lambda _{a})\ln(w_{j}-z_{a}) -\sum _{s< j}(\alpha_{i_{j}},\alpha_{i_{s}})\ln( w_{j}-w_{s})-\sum _{a< b}(\lambda_{a},\lambda_{b})\ln(z_{a}-z_{b})\\-c(\sum_{j}w_{j}-\frac{1}{2}\sum_{a}\parallel \lambda_{a}\parallel z_{a}).
\end{multline}
When $g=sl(2,\mathbb{C})$, there is only one simple root $\alpha$ satisfying $(\alpha,\alpha)=2$ and $(\alpha,\lambda_{a})=\lambda_{a}$. The generalized Bethe equation and Yang-Yang function degenerates to the Bethe equation \eqref{eq:a} and Yang-Yang function \eqref{eq:b} respectively and the highest weight $\lambda_{a}$ degenerates to the charge of singularity $z_{a}$. Now we use them as our start point to derive the braiding matrix of the thimble of the generalized Yang-Yang function.

\subsection{Definition and example of thimble}

We first introduce the concept of thimble. Before that we state two useful facts about the Morse function as the real part of a holomorphic function on an Hermitian manifold.
\begin{proposition}
\label{prop:a}
For a holomorphic function on an Hermitian manifold $M$ $dim_{\mathbb{R}}M=2d$, if its real part $h$ is a Morse function on $M$ (i.e. the Hessian matrix of $h$ is non-degenerate and its critical points are isolated), then (i) the gradient flow of $h$ keeps the imaginary part invariant; (ii) the index of each critical point of $h$ is $d$.
\end{proposition}
Proof:(i) Assuming that the holomorphic function is $f(x,y)=h(x,y)+ig(x,y)$ and the Hermitian metric on $M$ is $ds^{2}=H_{i\bar{j}}dz^{i}d\bar{z}^{j}$, then the gradient flow equation of $h$ is
\begin{equation}
\dfrac{dz^{i}}{dt}=-H^{i\bar{j}}\frac{\partial h}{\partial \bar{z}^{j}},
\end{equation}
\begin{equation}
\dfrac{d\bar{z}^{i}}{dt}=-\dfrac{\partial h}{\partial z^{j}}H^{j\bar{i}}.
\end{equation}
Thus \begin{equation}
\dfrac{dg}{dt}=\dfrac{1}{2i}\dfrac{d(f-\bar{f})}{dt}=\dfrac{1}{2i}(f_{z^{i}}\dfrac{dz^{i}}{dt}-\bar{f}_{\bar{z}^{i}}\dfrac{d\bar{z}^{i}}{dt})=\dfrac{1}{2i}(-f_{z^{i}}H^{i\bar{j}}\dfrac{\partial h}{\partial \bar{z}^{j}}+\dfrac{\partial h}{\partial z^{j}}H^{j \bar{i}}\bar{f}_{\bar{z}^{i}})=0.
\end{equation}
The last equality comes from the definition of the Hermitian metric.

(ii) From Cauchy-Rieaman equation, the Hessian matrix of $h$ is
\begin{equation}
\newline \left(
                                                                 \begin{array}{cc}
                                                                   h_{xx} & h_{xy} \\
                                                                   h_{xy} & h_{yy} \\
                                                                 \end{array}
                                                               \right)=\left(
                                                                         \begin{array}{cc}
                                                                           g_{xy} & g_{yy} \\
                                                                           -g_{xx} & -g_{xy} \\
                                                                         \end{array}
                                                                       \right),
\end{equation} thus $h_{xx}=-h_{yy}$. $h$ is a Morse function implies that the Hessian matrix is non-degenerate \begin{equation}-|h_{xx}|^{2}-|h_{xy}|^{2}\neq0,\end{equation} i.e. \begin{equation}|h_{xx}|\neq0\quad\mathrm{or}\quad |h_{xy}|\neq0.\end{equation}
\begin{description}
\item If $|h_{xx}|\neq0$, \begin{equation}
\left(
\begin{array}{cc}
                                                                   h_{xx} & h_{xy} \\
                                                                   h_{xy} & -h_{xx} \\
                                                                 \end{array}
                                                               \right)=\left(
                                                                         \begin{array}{cc}
                                                                           I & h_{xx}^{-1}h_{xy} \\
                                                                           0 & I \\
                                                                         \end{array}
                                                                       \right)^{t}
                                                               \left(
                                                                 \begin{array}{cc}
                                                                   h_{xx} & 0 \\
                                                                   0 & -h_{xy}h_{xx}^{-1}h_{xy}-h_{xx} \\
                                                                 \end{array}
                                                               \right)\left(
                                                                         \begin{array}{cc}
                                                                           I & h_{xx}^{-1}h_{xy} \\
                                                                           0 & I \\
                                                                         \end{array}
                                                                       \right).\end{equation}
                                                                        we consider the determinant of the matrix $(s\cdot h_{xy})h_{xx}^{-1}(s\cdot h_{xy})+h_{xx}$ with $s\in [0,1]$. $-|h_{xx}|^{2}-|h_{xy}|^{2}\neq0$ implies that
                                                                        \begin{equation}
                                                                        |(s\cdot h_{xy})h_{xx}^{-1}(s\cdot h_{xy})+h_{xx}|\neq0\end{equation}
                                                                         for any $s\in [0,1]$. $|(s\cdot h_{xy})h_{xx}^{-1}(s\cdot h_{xy})+h_{xx}|$ is a continuous function of $s$. Therefore, every eigenvalue of $(s\cdot h_{xy})h_{xx}^{-1}(s\cdot h_{xy})+h_{xx}$ keeps its sign invariant with $s$ varying from $0$ to $1$. $h_{xy}h_{xx}^{-1}h_{xy}+h_{xx}$ and $h_{xx}$ have the same index of inertia. $-h_{xy}h_{xx}^{-1}h_{xy}-h_{xx}$ and $h_{xx}$ have the opposite index of inertia. Thus the negative index of inertia of \begin{equation}\left(
                                                                 \begin{array}{cc}
                                                                   h_{xx} & h_{xy} \\
                                                                   h_{xy} & -h_{xx} \\
                                                                 \end{array}
                                                               \right)$$ is $d$. The index of the critical point is $d$.

                                                               \item If $|h_{xy}|\neq0$, $$2\left(
                                                                 \begin{array}{cc}
                                                                   h_{xx} & h_{xy} \\
                                                                   h_{xy} & -h_{xx} \\
                                                                 \end{array}
                                                               \right)=\left(
                                                                         \begin{array}{cc}
                                                                           I & I \\
                                                                           I & -I \\
                                                                         \end{array}
                                                                       \right)^{t}
                                                               \left(
                                                                 \begin{array}{cc}
                                                                   h_{xy} & h_{xx} \\
                                                                   h_{xx} & -h_{xy}  \\
                                                                 \end{array}
                                                               \right)\left(
                                                                         \begin{array}{cc}
                                                                           I & I \\
                                                                           I & -I \\
                                                                         \end{array}
                                                                       \right).
                                                                       \end{equation}
                                                                        With the same method above, we can prove that the index of the critical point is $d$.
                                                                       \end{description} This completes the proof.

Now we give a definition of thimble:
\begin{definition}
Assuming that $h$ as the real part of a holomorphic function on an Hermitian manifold $M$ is a Morse function on $M$. The cycle is called  a $thimble$ associated to $I$, denoted by $\mathcal {J}$, if all the points of it can be reached by the gradient flow of $h$ started from a critical point $I$.

\end{definition}

\begin{example}
A simple example of thimble comes from~\cite{h}:
$f(x)=i\lambda ( \dfrac {x^{3}} {3}-x)$ is called Airy function, where $\lambda$ is a constant in $\mathbb{C}$ and $x\in \mathbb{C}$.
The critical points of the holomorphic function $f$ are $x=\pm 1$, denoted by $P_{\pm}=\pm 1$. Assuming that the imaginary part of $\lambda$ is positive: $\lambda=a+bi$, where $b>0$, then $f(x)=-b(\dfrac {x^{3}} {3}-x)+ia(\dfrac {x^{3}} {3}-x)$. $\mathrm{Im} f(P_{+})=-\frac{2}{3}a$ and $\mathrm{Im} f(P_{-})=\frac{2}{3}a$. $\mathrm{Im} f(P_{+})=\mathrm{Im} f(P_{-})$ if and only if $a=0$. Thus, from the Proposition \ref{prop:a} above, there is a gradient flow connecting $P_{+}$ with $P_{-}$ if and only if $a=0$. When $a=0$, $f(x)=-b(\dfrac {x^{3}} {3}-x)$. The gradient flow connecting $P_{+}$ with $P_{-}$ is on the real axis of $x$ plane and the imaginary part of $f$ is zero along this flow. When $a\neq0$, there is no gradient flow connecting $P_{+}$ with $P_{-}$. As is shown in the Figure \ref{fig:i}, the picture (a), (b) and (c) describes the gradient flow started from $P_{+}$ and $P_{-}$ with $a=1$, $a=0$ and $a=-1$ respectively.

\begin{figure}[tbp]
\centering
\includegraphics[width=.6\textwidth,trim=0 200 0 200,clip]{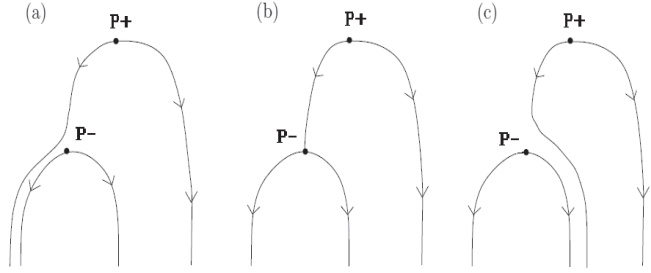}
\caption{\label{fig:i} Wall-crossing.}
\end{figure}

 If $a=1$ is continuously changed into $a=-1$ on the $\lambda$ plane with $b>0$, the thimble $\mathcal{J}_{+}$ associated to $P_{+}$ and the thimble $\mathcal{J}_{-}$ associated to $P_{-}$ will be transformed into $\mathcal{J}_{+}^{'}$ and  $\mathcal{J}_{-}^{'}$:\begin{equation}\left(
                                                                               \begin{array}{c}
                                                                                 \mathcal{J}_{+}^{'} \\
                                                                                 \mathcal{J}_{-}^{'} \\
                                                                               \end{array}
                                                                             \right)=\left(
                                                                                       \begin{array}{cc}
                                                                                         1 & \pm 1 \\
                                                                                         0 & 1 \\
                                                                                       \end{array}
                                                                                     \right)\left(
                                                                               \begin{array}{c}
                                                                                 \mathcal{J}_{+} \\
                                                                                 \mathcal{J}_{-} \\
                                                                               \end{array}
                                                                             \right).
                                                                             \end{equation}
$\mathcal{J}_{-}$ is invariant, but $\mathcal{J}_{+}$ will produce an additional term. This phenomena also appears when we continuously change $a=1$ into $a=-1$ with $b<0$.
When $a=0$, two rays $b>0$ and $b<0$ on the $\lambda$ plane are called Stokes rays (Stokes walls) by physicists. Passing through the Stokes ray is called wall-crossing. When wall-crossing happens, the thimble will produce an additional term.
\end{example}

\subsection{Braiding of thimble}

Now we focus on the braiding of the thimble coming from the real part of the Yang-Yang function
 $\mathscr{W}$. After the projection on the plane, knots can be decomposed as the contraction of the
 interaction, creation and annihilation operator (on page 117-118 of~\cite{c}). So it is enough to consider
 just two vertex operators or two singularities $z_{1}$ and $z_{2}$. Here we only consider the case
 $\lambda_{1}=\lambda_{2}$. There is no loss of generality in assuming that $z_{1}$ and $z_{2}$ have
 the same real part and $\mathrm{Im} z_{1}>\mathrm{Im} z_{2}$. Then we rotate the $z_{1}$ and $z_{2}$ clockwise by $\pi$
 around the middle point. In the case of no wall-crossing, $z_{1}$ and $z_{2}$ will change their position
 and the multiple valued function $e^{-\frac{\mathscr{W}}{k+h^{\vee }}}$ on the integration cycle $\Gamma$
 will be multiple of some power of $q$, where $q=e^{\frac{2\pi i}{k+h^{\vee}}}$. The original integration
 $\int _{\Gamma }e^{-\frac{\mathscr{W}}{k+h^{\vee }}}\prod _{j}dw_{j}$ becomes
 $q^{\theta(B,\lambda_{a},\alpha_{i_{j}})}\int _{\Gamma' }e^{-\frac{\mathscr{W}}{k+h^{\vee }}}\prod _{j}dw_{j}$, where $\Gamma'$ is the new integration cycle after braiding. The braiding transformation will be denoted by $\mathcal {B}$.
\begin{equation}
\mathcal {B}(\Gamma )=q^{\theta(B,\lambda_{a},\alpha_{i_{j}})}\Gamma',
\end{equation}
where the phase factor $\theta$ is a real number relative to the braiding and the weight of the highest weight representation $V_{\lambda_{a}}$ of the Lie algebra $g$. When the integration cycle $\Gamma$ is chosen to be a thimble, the phase factor can be computed easily, as is shown in the following lemma.

\begin{lemma}
\label{lemma:b}
If $\Gamma$ is a thimble associated to the real part of Yang-Yang function$\mathscr{W}(z_{1},z_{2},w_{j})$, then the phase factor of the integral $\int _{\Gamma }e^{-\frac{\mathscr{W}}{k+h^{\vee }}}\prod _{j}dw_{j}$ coming from the braiding without wall-crossing is equal to the phase factor of $e^{-\frac{\mathscr{W}_c}{k+h^{\vee }}}$ under braiding, where $\mathscr{W}_c$ is the value of the Yang-Yang function $\mathscr{W}(z_{1},z_{2},w)$ at the critical point $w_c$, i.e. $\mathscr{W}_c=\mathscr{W}(z_{1},z_{2},w_{c})$.
\end{lemma}
Proof: We assume that $W$ and $I$ is the real part and imaginary part of $\mathscr{W}$, i.e.$\mathscr{W}=W+iI$. Because the thimble is defined from the gradient flow of the real function $W$, from the first conclusion of the Proposition \ref{prop:a}, $I$ is a constant on the thimble.
From $\int _{\Gamma }e^{-\frac{\mathscr{W}}{k+h^{\vee }}}\prod _{j}dw_{j}=e^{-\frac{iI}{k+h^{\vee }}}\int _{\Gamma }e^{-\frac{W}{k+h^{\vee }}}\prod _{j}dw_{j}$, we see that the phase factor under braiding is coming from the braiding of $e^{-\frac{iI}{k+h^{\vee }}}$. And the function $e^{-\frac{iI}{k+h^{\vee }}}$ has the same phase factor with the function $e^{-\frac{\mathscr{W}_c}{k+h^{\vee }}}$ under the braiding. This concludes the proof.

First, we use Lemma \ref{lemma:b} to compute the phase factor of the braiding transformation without symmetry breaking.
\begin{equation}
e^{-\frac{\mathscr{W}}{k+h^{\vee }}}=\prod _{j}(w_{j}-z_{1})^{-\frac{(\alpha_{i_{j}},\lambda_{1})}{k+h^{\vee }}}(w_{j}-z_{2})^{-\frac{(\alpha_{i_{j}},\lambda_{2})}{k+h^{\vee }}}(z_{1}-z_{2})^{\frac{(\lambda_{1},\lambda_{2})}{k+h^{\vee }}}\prod_{j<s}(w_{j}-w_{s})^{\frac{(\alpha_{i_{j}},\alpha_{i_{s}})}{k+h^{\vee }}}.
\end{equation}
The Bethe equation is
\begin{equation}
\frac{(\alpha_{i_{j}},\lambda _{1})}{w_{j}-z_{1}}+\frac{(\alpha_{i_{j}},\lambda _{2})}{w_{j}-z_{2}}=0, j=1,2,...p,
\end{equation} with only one solution:
\begin{equation}
w_{j}=\frac{(\alpha_{i_{j}},\lambda _{2})z_{1}+(\alpha_{i_{j}},\lambda _{1})z_{2}}{(\alpha_{i_{j}},\lambda _{1})+(\alpha_{i_{j}},\lambda _{2})},j=1,2,...p.
\end{equation}
The corresponding thimble $\mathcal{J}_{p}$ connects $z_{1}$ with $z_{2}$
and passes through the critical point $w_{j}$. Thus both of factors $(w_{j}-z_{1})^{-\frac{(\alpha_{i_{j}},\lambda_{1})}{k+h^{\vee }}}$ and $(w_{j}-z_{2})^{-\frac{(\alpha_{i_{j}},\lambda_{2})}{k+h^{\vee }}}$ will have a contribution to the total phase factor. The thimble $\mathcal{J}_{p}$ after the braiding is still $\mathcal{J}_{p}$ but multiplied with some power of $q$. Using Lemma \ref{lemma:b}, we get the result as following:
   \begin{equation}\mathcal {B}e^{-\frac{W_{c}}{k+h^{\vee }}}=q^{-\frac{1}{2}[(\lambda_{1},\lambda_{2})+\sum_{j<s}(\alpha_{i_{j}},\alpha_{i_{s}})-\sum_{j,a}(\alpha_{i_{j}},\lambda_{a})]}e^{-\frac{W_{c}}{k+h^{\vee }}},\end{equation}
   thus
   \begin{equation}
   \label{eq:d}
   \mathcal{B}\mathcal{J}_{p}=(-1)^{p}q^{-\frac{1}{2}[(\lambda_{1},\lambda_{2})+\sum_{j<s}(\alpha_{i_{j}},\alpha_{i_{s}})-\sum_{j,a}(\alpha_{i_{j}},\lambda_{a})]}\mathcal{J}_{p}.
   \end{equation}
The $(-1)^{p}$ comes from the fact that the braiding changes the direction of each dimension of the thimble into the opposite direction and the thimble $\mathcal{J}_{p}$ is $p$ dimensional.

\begin{figure}[tbp]
\centering
\includegraphics[width=.6\textwidth,trim=0 200 0 200,clip]{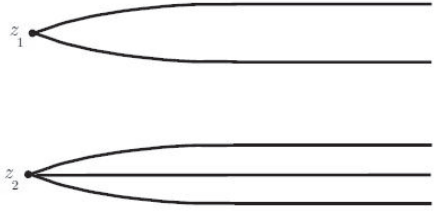}
\caption{\label{fig:j} Thimble $\mathcal{J}_{2,3}$.}
\end{figure}

Next, we consider the thimble of the real part of the generalized Yang-Yang function with symmetry breaking ($c>0$):
\begin{multline}
\mathscr{W}(w_{j},z_{1},z_{2})
=\sum _{j}(\alpha_{i_{j}},\lambda _{1})\ln(w_{j}-z_{1})+ \sum _{j}(\alpha_{i_{j}},\lambda _{2})\ln( w_{j}-z_{2})-\sum _{j<s}(\alpha_{i_{j}},\alpha_{i_{s}})\ln(w_{j}-w_{s})\\-(\lambda_{1},\lambda_{2})ln(z_{1}-z_{2})-c(\sum_{j}w_{j}-\frac{1}{2}\sum_{a}\parallel \lambda_{a}\parallel z_{a}).
\end{multline}
The Bethe equation is
\begin{equation}
\frac{(\alpha_{i_{j}},\lambda _{1})}{w_{j}-z_{1}}+\frac{(\alpha_{i_{j}},\lambda _{2})}{w_{j}-z_{2}}=\sum_{s\neq j}\frac{(\alpha_{i_{j}},\alpha _{i_{s}})}{w_{j}-w_{s}}+c, j=1,2,...p.
\end{equation}
When $c\rightarrow +\infty$, $w_{j}$ tends either to $z_{1}$ or $z_{2}$. It has several different solutions:
for any $s\in \mathbb{Z},0\leq s\leq p$,
\begin{equation}
w_{j}=\left\{
        \begin{array}{ll}
          z_{1}+o (\frac{1}{c}), & \hbox{$1\leq j\leq s$;} \\
         z_{2}+o(\frac{1}{c}), & \hbox{$(s+1)\leq j\leq p$,}
        \end{array}
      \right.
\end{equation}when $c$ is large enough.
To find solutions as vectors in the representation space $V_{\lambda}$ of the highest weight
representation of the Lie algebra $g$, we consider special solutions with the condition
$0\leq p\leq2(m-1)$ and $0\leq s\leq m-1$, where $m=\mathrm{dim}V_{\lambda}$.
Also $\lambda_{1}-\alpha_{i_{1}}-\alpha_{i_{2}}-...-\alpha_{i_{s}}$ and $\lambda_{2}-\alpha_{i_{s+1}}-\alpha_{i_{s+2}}-...-\alpha_{i_{p}}$
should be the weights of the representation $V_{\lambda}$. We denote the thimble of this kind as
$\mathcal{J}_{s,p-s}$. It is a $p$ dimensional sub-manifold in a $2p$
dimensional manifold $\mathbb{C}^{p}$.
Therefore, when $c\rightarrow +\infty$, the gradient flow associated to $w_{j}$ is from
$z_{1}$ to infinity for $1\leq j\leq s$ and from $z_{2}$ to infinity for $(s+1)\leq j\leq p$.
The thimble $\mathcal{J}_{s,p-s}$ is the Cartesian product of such $p$ one dimensional
 manifolds. For example, $\mathcal{J}_{2,3}$ is a $5$ dimensional manifold as the Cartesian product of $5$
 gradient flows, $2$ from $z_{1}$ and $3$ from $z_{2}$, as is shown in the Figure \ref{fig:j}.

Now combining all thimbles together with respect to $p$ from $0$ to $2(m-1)$, we have totally
$m^{2}$ different thimbles. Each thimble $\mathcal{J}_{s,p-s}$ corresponds to a weight vector
with weight $\lambda_{1}-\sum_{1\leq j\leq s}\alpha_{i_{j}},\lambda_{2}-\sum_{s+1\leq j\leq p}\alpha_{i_{j}}$
in the representation space $V_{\lambda_{1}}\otimes V_{\lambda_{2}}$. Therefore, all the $m^{2}$ thimbles
as special solutions of the Bethe equation in the symmetry breaking case naturally form a set of bases
of the representation space $V_{\lambda_{1}}\otimes V_{\lambda_{2}}$.
With these bases the representation space $V_{\lambda_{1}}\otimes V_{\lambda_{2}}$ is naturally decomposed into
the direct summation of the subspace with respect to the dimension $p$ of the thimbles:
\begin{equation}
V_{\lambda_{1}}\otimes V_{\lambda_{2}}=\oplus_{p=0}^{2m-2}V_{p},
\end{equation}
where $V_{p}\triangleq \{v_{s}\otimes v_{p-s}, 0\leq s\leq p \mid v_{s}\otimes v_{p-s}\in V_{\lambda_{1}}\otimes V_{\lambda_{2}}, v_{s}\otimes v_{p-s} \hbox {is a vector with weight} (\lambda-\alpha_{i_{1}}-...-\alpha_{i_{s}},\lambda-\alpha_{i_{s+1}}-...-\alpha_{p})  \}\subseteq V_{\lambda_{1}}\otimes V_{\lambda_{2}}$.
Clearly the braiding does not change the dimension of the thimble, thus each $V_{p}$ is an invariant subspace for the  braiding operator $\mathcal {B}$.

 For thimble $\mathcal{J}_{p,0}$, the clockwise braiding will simply change this thimble into the thimble $\mathcal{J}_{0,p}$ up to an phase factor.
 Because $w_{j}$ is near $z_{1}$, by Lemma \ref{lemma:b}, the factor $(w_{j}-z_{1})^{-\frac{(\alpha_{i_{j}},\lambda_{1})}{k+h^{\vee }}}$ and $(w_{j}-w_{s})^{\frac{(\alpha_{i_{j}},\alpha_{i_{s}})}{k+h^{\vee }}}$ do not have any contribution to the phase factor, but the factor $(w_{j}-z_{2})^{-\frac{(\alpha_{i_{j}},\lambda_{2})}{k+h^{\vee }}}$ does. Thus \begin{equation}\mathcal {B}\mathcal{J}_{p,0}=q^{-\frac{1}{2}(\lambda_{1}-\sum_{j}\alpha_{i_{j}},\lambda_{2})}\mathcal{J}_{0,p}.
\end{equation}
It should be noticed that $\lambda_{1}-\sum_{j}\alpha_{i_{j}}$ and $\lambda_{2}$ are the $(p+1)$th and the first weight of the representation $V_{\lambda_{1}}$ and $V_{\lambda_{2}}$ respectively. When thimble is connecting $z_{2}$ to the infinity, for example $\mathcal{J}_{0,1}$, the braiding will cause additional wall-crossing terms as indicated before. However, the non-wall-crossing term still can be determined by Lemma \ref{lemma:b}. As the similar discussion above, for thimble $\mathcal{J}_{s,p-s}(p-s\neq0)$,
\begin{equation}
\label{eq:e}
\mathcal {B}\mathcal{J}_{s,p-s}=q^{-\frac{1}{2}(\lambda_{1}-\sum_{1\leq j\leq s}\alpha_{i_{j}},\lambda_{2}-\sum_{s+1\leq j\leq p}\alpha_{i_{j}})}\mathcal{J}_{p-s ,s}+ w.c.t.,
\end{equation}
where w.c.t are undetermined wall-crossing terms.

Two properties of wall-crossing should be noticed: First,
The braiding transformation keeps invariant the total types and numbers of the simple roots associated to
the thimble it acts on. If the total types and numbers of the simple roots of two thimbles are different,
the Yang-Yang functions of them are two different functions. From the definition of the braiding of the
thimble, the braiding transformation only acts on the thimbles from one holomorphic function. Thus this
property is natural. Wall-crossing in the braiding transformation does not create or annihilate any simple roots, but only transfers
 them from one location to another.
We call this property the conservation law of wall-crossing; Second, the transfer of simple roots in the wall-crossing can
only be from $z_{2}$ to $z_{1}$. The gradient flows in symmetry breaking case are from $z_{1}$ and
$z_{2}$ to the infinity in the positive direction of the real axis in the $w$ plane. In our assumption, $z_{1}$ and $z_{2}$ have
 the same real part and $\mathrm{Im} z_{1}>\mathrm{Im}z_{2}$. Therefore, in the clockwise braiding, the wall-crossing appears when there is a gradient flow from $z_{2}$ to the infinity passing through $z_{1}$.
Thus the only possible transfer of simple roots is from $z_{2}$
 to $z_{1}$. These two properties tell us the braiding matrix is a diagonal
partitioned matrix and each block in the diagonal is a triangular matrix with respect to each $p$.
 Thus, they are actually sub-representations of the braiding. For example, the block of
 $p\leq m-1$ is
\begin{multline}
\mathcal {B}\left(
              \begin{array}{c}
                \mathcal{J}_{p,0} \\
                \mathcal{J}_{p-1,1} \\
                : \\
                \mathcal{J}_{0,p} \\
              \end{array}
            \right)\\=\left(
                      \begin{array}{cccc}
                          &   &   & q^{-\frac{1}{2}(\lambda_{1}-\sum_{1\le j\le p }\alpha_{i_{j}},\lambda_{2})} \\
                          &   & q^{-\frac{1}{2}(\lambda_{1}-\sum_{1\le j\le p-1 }\alpha_{i_{j}},\lambda_{2}-\alpha_{i_{p}})} & * \\
                          & : & * & * \\
                        q^{-\frac{1}{2}(\lambda_{1},\lambda_{2}-\sum_{1\le j\le p }\alpha_{i_{j}})} & * & * & * \\
                      \end{array}
                    \right)\left(
              \begin{array}{c}
                \mathcal{J}_{p,0} \\
                \mathcal{J}_{p-1,1} \\
                : \\
                \mathcal{J}_{0,p} \\
              \end{array}
            \right).
\end{multline}
Every skew diagonal element of the diagonal block of the braiding matrix is derived from Lemma \ref{lemma:b}. And its phase factor is coming from the inner product of two weights of representation $V_{\lambda_{1}}$ and $V_{\lambda_{2}}$.

When skew diagonal elements in the triangular matrix are known, the remaining problem is to find the wall-crossing term. In~\cite{e}, integration cycles are used to compute the braiding
transformation, then they are transformed into bases of thimbles to get the braiding matrix needed. Now we use the same method to derive the braiding matrix of the fundamental representation of $sl(n+1,\mathbb{C})$. Here we use Dynkin label, $\lambda_{1}=\lambda_{2}=\lambda=(1,0,...,0)$.
$sl(n+1,\mathbb{C})$ has $n$ simple roots $\alpha_{i}, i=1,2,...,n$.
The weights of the fundamental representation are $\lambda, (\lambda-\alpha_{1}), (\lambda-\alpha_{1}-\alpha_{2}),(\lambda-\alpha_{1}-\alpha_{2}-\alpha_{3}),...,(\lambda-\alpha_{1}-\alpha_{2}-...-\alpha_{n}).$
We denote them as $\lambda^{0},\lambda^{1},...,\lambda^{n}$.

\begin{lemma}
\label{lemma:c}
\begin{equation}
(\lambda^{s},\lambda^{t})=\left\{
                            \begin{array}{ll}
                              \frac{n}{n+1}, & \hbox{$s=t$;} \\
                              -\frac{1}{n+1}, & \hbox{$s\neq t$.}
                            \end{array}
                          \right.
\end{equation}
\end{lemma}
Proof: By straightforward calculation.

\begin{lemma}
\begin{equation}
(\lambda^{s},\lambda^{t})=(\lambda^{n-s},\lambda^{n-t}).
\end{equation}
\end{lemma}
Proof: $n-s=n-t$ if and only if $s=t$. From Lemma \ref{lemma:c}, the proof is straightforward.

This duality property for the fundamental representation of $sl(n+1,\mathbb{C})$ implies that blocks of $p=i$ and $p=2n-i$ in the braiding matrix are same.

\begin{theorem}
\label{Thm:braiding}

For general $p$ of $0\leq p\leq 2n$, $0\leq m\leq n$ and $0\leq p-m\leq n$:
\begin{equation}
\mathcal{B}\mathcal{J}_{m, p-m}=\left\{
                                             \begin{array}{ll}
                                               q^{- \frac{n}{2(n+1)}}\mathcal{J}_{m, p-m}, & \hbox{$m= p-m$;} \\
                                               q^{ \frac{1}{2(n+1)}}\mathcal{J}_{ p-m,m}, & \hbox{$m > p-m$;} \\
                                               q^{ \frac{1}{2(n+1)}}\mathcal{J}_{ p-m,m}+(q^{-\frac{n}{2(n+1)}}-q^{\frac{n+2}{2(n+1)}})\mathcal{J}_{m, p-m}, & \hbox{$m< p-m$.}
                                             \end{array}
                                           \right.
\end{equation}

\end{theorem}
Proof:\begin{itemize}
        \item $m= p-m$: There is no wall-crossing. From Lemma \ref{lemma:b} and Lemma \ref{lemma:c}, $$\mathcal{B}\mathcal{J}_{m, p-m}=q^{-\frac{1}{2}(\lambda^{m},\lambda^{m})}\mathcal{J}_{m, p-m}=q^{- \frac{n}{2(n+1)}}\mathcal{J}_{m, p-m};$$
        \item  $m >  p-m$: Also there is no wall-crossing. From Lemma \ref{lemma:b} and Lemma \ref{lemma:c},
            $$\mathcal{B}\mathcal{J}_{m, p-m}=q^{-\frac{1}{2}(\lambda^{m},\lambda^{ p-m})}\mathcal{J}_{m, p-m}=q^{ \frac{1}{2(n+1)}}\mathcal{J}_{m, p-m};$$
        \item  $m <  p-m$: From the conservation law of wall-crossing, there will be one wall-crossing term of $\mathcal{J}_{m, p-m}$. From formula (\ref{eq:e}) and Lemma \ref{lemma:c}, $$\mathcal{B}\mathcal{J}_{m, p-m}=q^{ -\frac{1}{2}(\lambda^{m},\lambda^{ p-m})}\mathcal{J}_{ p-m,m}+d\mathcal{J}_{m, p-m}=q^{ \frac{1}{2(n+1)}}\mathcal{J}_{ p-m,m}+d\mathcal{J}_{m, p-m},$$ where $d$ is a constant to be determined. The transformation of $\mathcal{J}_{m, p-m}$ and $\mathcal{J}_{ p-m,m}$ forms into a matrix:
          $$\mathcal{B}\left(
                         \begin{array}{c}
                           \mathcal{J}_{ p-m,m} \\
                           \mathcal{J}_{m, p-m} \\
                         \end{array}
                       \right)=\left(
                                 \begin{array}{cc}
                                   0 & q^{ \frac{1}{2(n+1)}} \\
                                   q^{ \frac{1}{2(n+1)}} & d \\
                                 \end{array}
                               \right)\left(
                                        \begin{array}{c}
                                          \mathcal{J}_{ p-m,m} \\
                                          \mathcal{J}_{m, p-m} \\
                                        \end{array}
                                      \right).
$$
To determine $d$, we derive the braiding matrix of cycles $C_{ p-m,m}$ and $C_{m, p-m}$. For convenience, we assume that $ p-m=m+l$. From the second property of wall-crossing, the only possible transfer of simple roots is from $z_{2}$
 to $z_{1}$ , so braiding of $C_{ p-m,m}$ is easy: $$\mathcal{B}C_{ p-m,m}=q^{-\frac{1}{2}(\lambda^{0},\lambda^{0})}C_{m, p-m}=q^{-\frac{n}{2(n+1)}}C_{m, p-m}.$$

\begin{figure}[tbp]
\centering
\includegraphics[width=.8\textwidth,trim=0 200 0 200,clip]{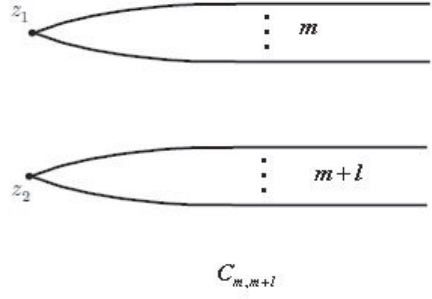}
\caption{\label{fig:Cmm+l}$C_{m,m+l}$ before braiding}
\end{figure}

\begin{figure}[tbp]
\centering
\includegraphics[width=.8\textwidth,trim=0 200 0 200,clip]{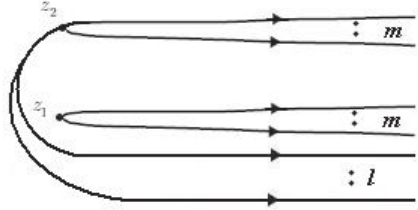}
\caption{\label{fig:WC1}$C_{m,m+l}$ after braiding}
\end{figure}

\begin{figure}[tbp]
\centering
\includegraphics[width=.7\textwidth,trim=0 200 0 200,clip]{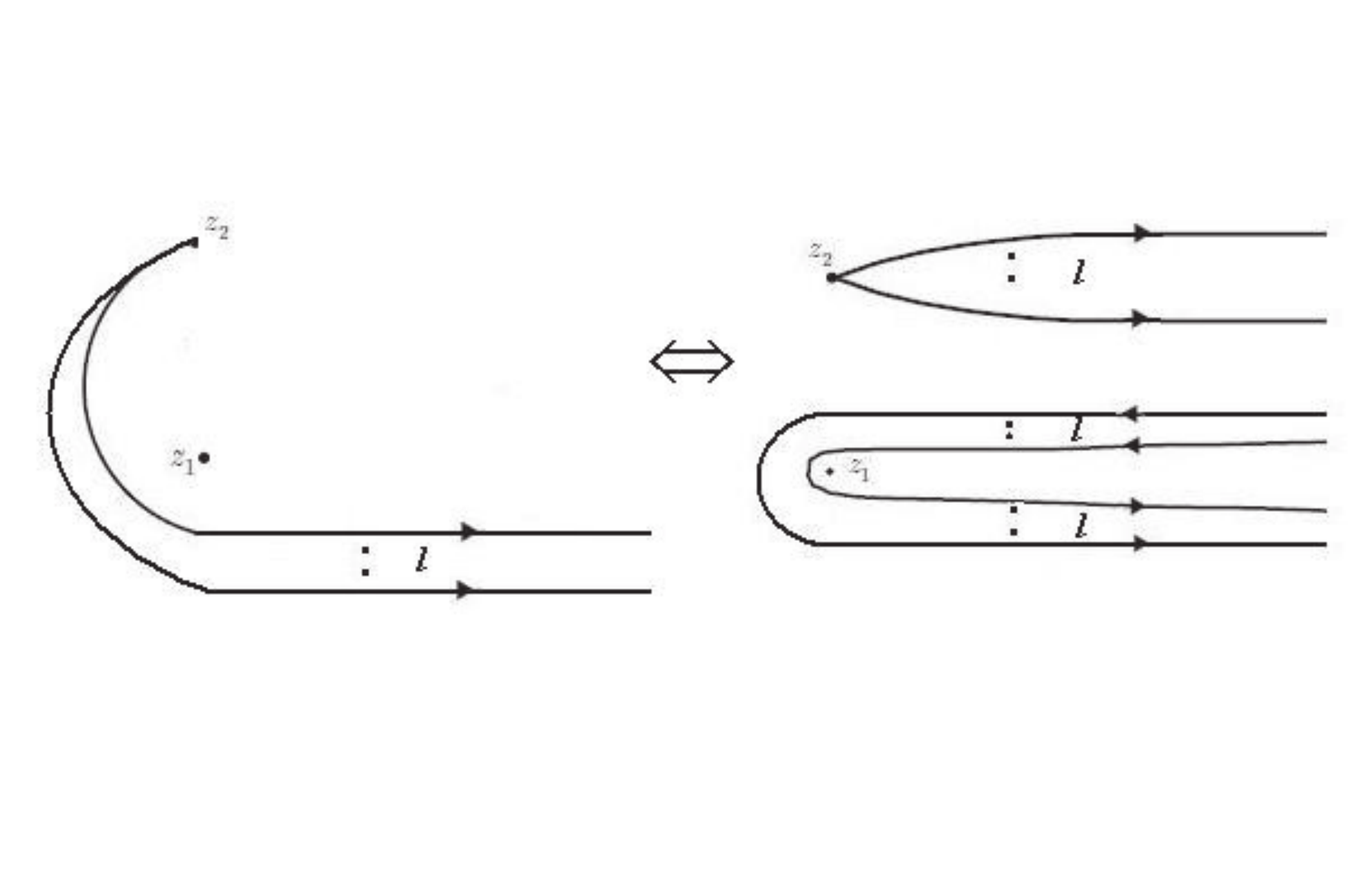}
\caption{\label{fig:WC2}Homology equivalence of wall-crossing part}
\end{figure}

The braiding of $C_{m, p-m}=C_{m,m+l}$ will cause wall-crossing. As is shown in the Figure \ref{fig:Cmm+l}, Figure \ref{fig:WC1} and Figure \ref{fig:WC2},  wall-crossing part is equivalent in homology to a zig-zag cycle, which is starting at $z_{2}$, heads directly to $Re z = \infty$ before doubling back around $z_{1}$ and returning to $Re z = \infty$. Thus, there are are three pieces in the wall-crossing part and two pieces near $z_{1}$ are different by a deck transformation $q\cdot$. Therefore, we have

\begin{equation}
\begin{split}
\mathcal{B}C_{m, p-m}
&=q\cdot q^{-\frac{1}{2}(\lambda^{0},\lambda^{0})}(C_{ p-m,m}-C_{m, p-m})+q^{-\frac{1}{2}(\lambda^{0},\lambda^{0})}C_{m, p-m}\\
&=q\cdot q^{ -\frac{n}{2(n+1)}}(C_{ p-m,m}-C_{m, p-m})+q^{ -\frac{n}{2(n+1)}}C_{m, p-m}\\
&=q^{ \frac{n+2}{2(n+1)}}C_{ p-m,m}+(q^{ -\frac{n}{2(n+1)}}-q^{ \frac{n+2}{2(n+1)}})C_{m, p-m},
\end{split}
\end{equation}
where $q\cdot$ is from deck transformation.

Thus,
$$\mathcal{B}\left(
                         \begin{array}{c}
                           C_{ p-m,m} \\
                           C_{m, p-m} \\
                         \end{array}
                       \right)=\left(
                                 \begin{array}{cc}
                                   0 & q^{-\frac{n}{2(n+1)}} \\
                                   q^{ \frac{n+2}{2(n+1)}} & q^{ -\frac{n}{2(n+1)}}-q^{ \frac{n+2}{2(n+1)}} \\
                                 \end{array}
                               \right)\left(
                                        \begin{array}{c}
                                          C_{ p-m,m} \\
                                          C_{m, p-m} \\
                                        \end{array}
                                      \right).
$$
$\{C_{ p-m,m}, C_{m, p-m}\}$ and $\{\mathcal{J}_{ p-m,m}, \mathcal{J}_{m, p-m}\}$ are two bases in the same vector space, braiding matrixes in these two bases are similar to each other. Thus $d=q^{ -\frac{n}{2(n+1)}}-q^{ \frac{n+2}{2(n+1)}}.$ This completes the proof.
      \end{itemize}

It should be noticed that from the braiding matrix above the irreducible representation of the  braiding on $V_{\lambda_{1}}\otimes V_{\lambda_{2}}$ is not larger than two dimension.

In the appendix, we have checked that the braiding matrices we derived satisfy Yang-Baxter equation. Thus, they give a representation of braid group $B_{2}$. The tensor product of braiding matrix and $m-2$ identities generates the representation of braid group $B_{m}$.

\subsection{Fusion matrix}

When the braiding matrix is known, the amplitudes for creation or annihilation of a pair of strands can be determined. As in the case of $g=sl(2,\mathbb{C})$ (see p75-76 in  ~\cite{e}), by the correspondence from vectors in $V_{\lambda_{1}}\otimes V_{\lambda_{2}}$ to thimbles, if we denote $+$ and $-$ to be the vector of weight $\lambda-\alpha$ and the vector of weight $\lambda$, then we can write all amplitudes for annihilation into a matrix:
                             $$\mathcal{M}=\left(
                                      \begin{array}{ccc}
                                           & + & - \\
                                          + & 0 & iq^{-\frac{1}{4}} \\
                                          - & -iq^{\frac{1}{4}} & 0 \\
                                        \end{array}
                                      \right),$$ where $i$ is a normalization constant so that two annihilation amplitudes between two complementary states are inverse to each other. All amplitudes for creation form another matrix, under
                                      the normalization above, they are same.
                                      \begin{figure}[tbp]
\centering
\includegraphics[width=.8\textwidth,trim=0 200 0 200,clip]{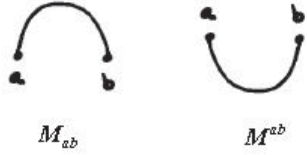}
\caption{\label{fig:AC} Annihilation $\mathcal{M}_{ab}$ and creation $\mathcal{M}^{ab}$}
\end{figure}
In general, we denote the amplitudes for annihilation and creation between two states $v_{a}$ and $v_{b}$ as $\mathcal{M}_{ab}$ and$\mathcal{M}^{ab}$ respectively ( see Figure \ref{fig:AC} ).
\begin{figure}[tbp]
\centering
\includegraphics[width=.8\textwidth,trim=0 200 0 200,clip]{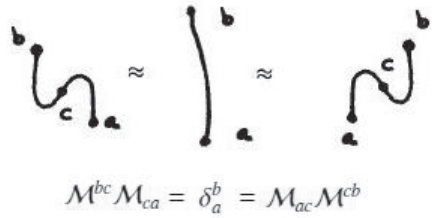}
\caption{\label{fig:TM}Invariance under the topological move}
\end{figure}

As is shown in the Figure \ref{fig:TM}, to be invariant under the topological move, they should be inverse to each other:
\begin{equation}
\mathcal{M}^{bc}\mathcal{M}_{ca}=\mathcal{M}_{ac}\mathcal{M}^{cb}=\delta ^{b}_{a}.
\end{equation}
We call the matrix of amplitudes for the annihilation of two strands fusion matrix $\mathcal{M}$, then the matrix of amplitudes for the creation of two strands is just its inverse $\mathcal{M}^{-1}$.
For the fundamental representation of $A_{n}$, we define two vectors $v_{s}$ of weight $\lambda-\alpha_{i_{1}}-...-\alpha_{i_{s}}$ and $v_{n-s}$ of weight $\lambda-\alpha_{j_{1}}-...-\alpha_{j_{n-s}}$ in the $V_{\lambda_{1}}$ and $V_{\lambda_{2}}$ to be complementary to each other. Since only two complementary states can fuse into a vacuum state,
fusion amplitudes are nonzero only between two complementary vectors.  Then fusion matrix can be written as
$\mathcal{M}:V_{\lambda}\longrightarrow V_{\lambda}$
            \begin{equation}\mathcal{M}\left(
                                                                            \begin{array}{c}
                                                                              v_{n}\\
                                                                              v_{n-1}\\
                                                                              ... \\
                                                                              v_{1} \\
                                                                              v_{0} \\
                                                                            \end{array}
                                                                          \right)
            =\left(
                                                 \begin{array}{ccccc}
                                                    &  &  &  &  m_{0}\\
                                                    & 0 &  & m_{1} &  \\
                                                    &  & ... &  &  \\
                                                    & m_{n-1} &  & 0 &  \\
                                                   m_{n} &  &  &  &  \\
                                                 \end{array}
                                               \right)\left(
                                                                            \begin{array}{c}
                                                                              v_{n}\\
                                                                              v_{n-1}\\
                                                                              ... \\
                                                                              v_{1} \\
                                                                              v_{0} \\
                                                                            \end{array}
                                                                          \right).
            \end{equation}

And its inverse is:
 \begin{equation}\mathcal{M}^{-1}\left(
                                                                            \begin{array}{c}
                                                                              v_{n}\\
                                                                              v_{n-1}\\
                                                                              ... \\
                                                                              v_{1} \\
                                                                              v_{0} \\
                                                                            \end{array}
                                                                          \right)
            =\left(
                                                 \begin{array}{ccccc}
                                                    &  &  &  &  m_{n}^{-1}\\
                                                    & 0 &  & m_{n-1}^{-1} &  \\
                                                    &  & ... &  &  \\
                                                    & m_{1}^{-1} &  & 0 &  \\
                                                   m_{0}^{-1} &  &  &  &  \\
                                                 \end{array}
                                               \right)\left(
                                                                            \begin{array}{c}
                                                                              v_{n}\\
                                                                              v_{n-1}\\
                                                                              ... \\
                                                                              v_{1} \\
                                                                              v_{0} \\
                                                                            \end{array}
                                                                          \right).
            \end{equation}

\begin{figure}[tbp]
\centering
\includegraphics[width=.8\textwidth,trim=0 200 0 200,clip]{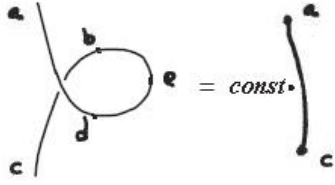}
\caption{\label{fig:Mcdt} Constraint for $\mathcal{M}$}
\end{figure}

To derive a knot invariant, as is shown in the Figure \ref{fig:Mcdt}, braiding and fusion matrix must satisfy the following condition:
\begin{equation}\label{equ:cond for M}
\sum_{b,d,e}\mathcal{B}_{cd}^{ab}\mathcal{M}_{be}\mathcal{M}^{de}=c\delta_{c}^{a}, \hbox{where c is a constant}.
\end{equation}

\begin{theorem}
For the braiding matrix in Theorem \ref{Thm:braiding}, the condition $$\sum_{b,d,e}\mathcal{B}_{cd}^{ab}\mathcal{M}_{be}\mathcal{M}^{de}=c\delta_{c}^{a}$$ is equivalent to  $$\frac{m_{a}}{m_{n-a}}=q^{a-\frac{n}{2}}\hbox{ or }\frac{m_{a}}{m_{n-a}}=-q^{a-\frac{n}{2}},\quad a=0,1,2,...,n.$$
\end{theorem}
Proof: The left hand side of (\ref{equ:cond for M}) is :
\begin{itemize}
  \item $a\neq c$ $$\hbox{L.H.S}=\sum_{e}\mathcal{B}_{cn-e}^{an-e}\mathcal{M}_{n-ee}\mathcal{M}^{n-ee}=0$$
  \item $a=c$
\begin{equation}
\begin{split}L.H.S&=\sum_{e\leq n-a}\mathcal{B}_{an-e}^{an-e}\mathcal{M}_{n-ee}\mathcal{M}^{n-ee}\\
                  &=q^{-\frac{n}{2(n+1)}}\mathcal{M}_{an-a}\mathcal{M}^{an-a}+(q^{-\frac{n}{2(n+1)}}-q^{\frac{n+2}{2(n+1)}})\sum_{e<n-a}\mathcal{M}_{n-ee}\mathcal{M}^{n-ee}
\end{split}
\end{equation}
\end{itemize}
Thus, the condition (\ref{equ:cond for M}) leads to a group of linear equations for $\frac{m_{e}}{m_{n-e}},\quad e=0,1,2,...,n$:
\begin{equation}
q^{-\frac{n}{2(n+1)}}\mathcal{M}_{an-a}\mathcal{M}^{an-a}+(q^{-\frac{n}{2(n+1)}}-q^{\frac{n+2}{2(n+1)}})\sum_{e<n-a}\mathcal{M}_{n-ee}\mathcal{M}^{n-ee}=c, \quad a=0,1,2,...,n.
\end{equation}
It has two solutions:
$$\frac{m_{a}}{m_{n-a}}=q^{a-\frac{n}{2}}\hbox{ or }\frac{m_{a}}{m_{n-a}}=-q^{a-\frac{n}{2}},\quad a=0,1,2,...,n.$$
This completes the proof.

Here we choose a special solution satisfying the following normalization condition:
\begin{equation}
m_{a}\cdot m_{n-a}=1,\quad a=0,1,2,...,n,
\end{equation}
or equivalently,
\begin{equation}\label{Norm for fusion}
\mathcal{M}_{ab}=\mathcal{M}^{ab},\quad a,b=0,1,2,...,n.
\end{equation}

\begin{theorem}
For the braiding matrix in Theorem \ref{Thm:braiding}, $M$ satisfies (\ref{equ:cond for M}) and normalization condition (\ref{Norm for fusion}) if and only if
\begin{equation}
\mathcal{M}_{ab}=\left\{
         \begin{array}{ll}
           0, & \hbox{$a+b\neq n$;} \\
           c_{0}q^{\frac{n-2a}{4}}, & \hbox{$a+b=n$,}
         \end{array}
       \right.
\end{equation} where $c_{0}=\pm1,\pm i$.
\end{theorem}
Proof: From the solutions of (\ref{equ:cond for M}), the normalization condition (\ref{Norm for fusion}) demands that $$(m_{a})^{2}=q^{a-\frac{n}{2}}\hbox{ or }(m_{a})^{2}=-q^{a-\frac{n}{2}},\quad a=0,1,2,...,n.$$ This leads to $$m_{a}=c_{0}q^{\frac{2a-n}{4}},\hbox{ i.e. }M_{ab}= c_{0}q^{\frac{n-2a}{4}}\delta_{a+b}^{n},\hbox{ where }c_{0}=\pm1,\pm i,\quad a=0,1,2,...,n.$$
This completes the proof.

For convenience, in the following we choose $c_{0}=1$, i.e.
\begin{equation}
\mathcal{M}_{ab}=\left\{
         \begin{array}{ll}
           0, & \hbox{$a+b\neq n$;} \\
           q^{\frac{n-2a}{4}}, & \hbox{$a+b=n$.}
         \end{array}
       \right.
\end{equation}

\subsection{Knots invariant from braiding and fusion matrix}
\begin{figure}[tbp]
\centering
\includegraphics[width=.9\textwidth,trim=0 200 0 200,clip]{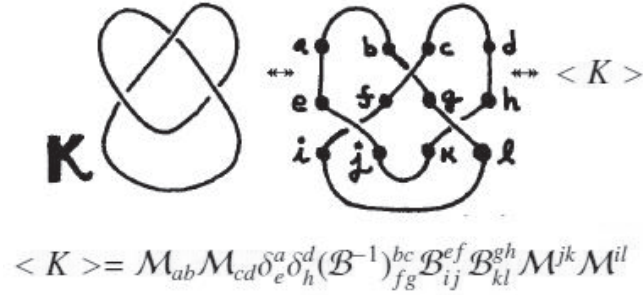}
\caption{\label{fig:knot}Decomposition of a link diagram}
\end{figure}

In \cite{c}, the method of quantum mechanics was used to study knots invariants. Every link in $\mathbb{R}^{3}$ can be projected on a plane as a link diagram. After that, it can be decomposed into the combination of  braiding, fusion and identity, as shown in the Figure \ref{fig:knot}. Then knots can be thought as a process with braiding, fusions or identities as its intermediate configurations. Thus knot invariant $<K>$ as an expectation of a quantum mechanics system is just the contraction of the braiding and fusion and their inverses appeared in the decomposition of the diagram of knot $K$. In the Figure \ref{fig:knot}, $$<K>=\mathcal{M}_{ab}\mathcal{M}_{cd}\delta_{e}^{a}\delta_{h}^{d}(\mathcal{B}^{-1})_{fg}^{bc}\mathcal{B}_{ij}^{ef}\mathcal{B}_{kl}^{gh}\mathcal{M}^{jk}\mathcal{M}^{il},$$
where we use Einstein notation for summation.

\begin{figure}[tbp]
\centering
\includegraphics[width=.8\textwidth,trim=0 200 0 200,clip]{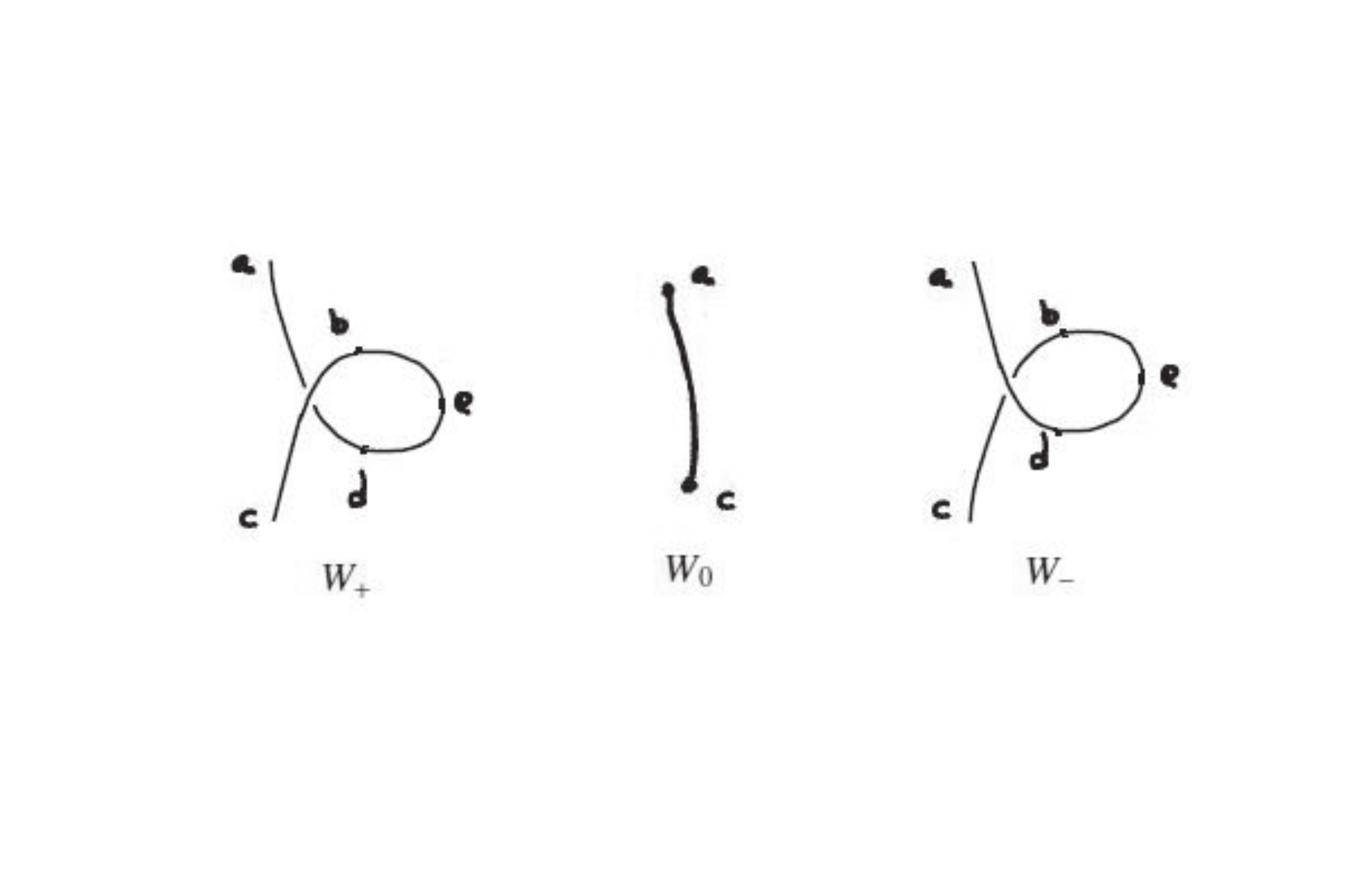}
\caption{\label{fig:UOT}}
\end{figure}
Here is an example:
\begin{example}
As is shown in the Figure \ref{fig:UOT}, $W_{+}$ and $W_{-}$ are twist and anti-twist respectively. Let $<K>$ be the invariant defined as an expectation of a quantum mechanics system associated to the knot $K$. Then

\begin{equation}
\begin{split}<W_{+}>&=\sum_{b,d,e}(\mathcal{B}^{-1})_{cd}^{ab}\mathcal{M}_{be}\mathcal{M}^{de}\\
                    &=\sum_{n-a\leq e}(\mathcal{B}^{-1})_{an-e}^{an-e}\mathcal{M}_{n-ee}\mathcal{M}^{n-ee}\delta_{c}^{a}\\
                    &=q^{\frac{n}{2(n+1)}}\mathcal{M}_{an-a}\mathcal{M}^{an-a}+(q^{\frac{n}{2(n+1)}}-q^{-\frac{n+2}{2(n+1)}})\sum_{e>n-a}\mathcal{M}_{n-ee}\mathcal{M}^{n-ee}\delta_{c}^{a}\\
                    &=q^{\frac{n}{2(n+1)}}\mathcal{M}_{an-a}\mathcal{M}^{an-a}+(q^{\frac{n}{2(n+1)}}-q^{-\frac{n+2}{2(n+1)}})q\cdot\mathcal{M}_{an-a}\mathcal{M}^{an-a}\frac{1-q^{a}}{1-q}\delta_{c}^{a}\\
                    &=q^{\frac{n}{2(n+1)}}\mathcal{M}_{an-a}\mathcal{M}^{an-a}(1+(1-q^{-1})\cdot q\frac{1-q^{a}}{1-q}\delta_{c}^{a})\\
                    &=q^{\frac{n}{2(n+1)}}\mathcal{M}_{an-a}\mathcal{M}^{an-a}\cdot q^{a}\delta_{c}^{a}\\
                    &=q^{\frac{n}{2(n+1)}}\mathcal{M}_{0n}\mathcal{M}^{0n}\delta_{c}^{a}\\
                    &=q^{\frac{n(n+2)}{2(n+1)}}\delta_{c}^{a}\\
                    &=q^{\frac{n(n+2)}{2(n+1)}}<W_{0}>
\end{split}
\end{equation}

Similarly,
\begin{equation}<W_{-}>=\sum_{b,d,e}\mathcal{B}_{cd}^{ab}\mathcal{M}_{be}\mathcal{M}^{de}=q^{-\frac{n(n+2)}{2(n+1)}}<W_{0}>.\end{equation}
This example shows that the knot invariant $<K>$ from our construction will produce an factor $a=q^{\frac{n(n+2)}{2(n+1)}}$ or $a^{-1}=q^{-\frac{n(n+2)}{2(n+1)}}$under Reidemeister move of type I. If we thicken a knot into a ribbon by natural framing, i.e. choosing a frame to be normal to the plane the knot projected on, then $<K>$ depends on its framing.
\end{example}

$<K>$ is called an isotopy invariant, if it is invariant under three kinds of Reidemeister moves. $<K>$ is called a regular isotopy invariant, if it is invariant under Reidemeister moves of type II and type III.

In the appendix, we prove that the expectations $<K>$ of a quantum mechanics system involving the braiding derived in Theorem \ref{Thm:braiding} and fusion satisfying (\ref{equ:cond for M}) is invariant under Reidemeister move of type II and type III. Thus, it is a regular isotopy invariant.

However, HOMFLY polynomial is an isotopy invariant. To derive HOMFLY polynomial, we define  $H(K)$ from $<K>$ as follows.

\begin{definition}
For an oriented knot K, define H(K) to be
\begin{equation}
H(K)=a^{-\omega(K)}<K>,
\end{equation}where $\omega(K)$ is the writhe of K.
\end{definition}
For example, $\omega(L_{+})=1$ and $\omega(L_{-})=-1$, then $$H(L_{+})=a^{-1}<L_{+}>=q^{-\frac{n(n+2)}{2(n+1)}}<L_{+}>\hbox{ and  } H(L_{-})=a<L_{+}>=q^{\frac{n(n+2)}{2(n+1)}}<L_{-}>.$$
It is easy to prove that $H(K)$ is independent of framing of $K$ (see Proposition 3.7 in ~\cite{c}). We prove that $K$ is a regular isotopy invariant in the appendix,  therefore $H(K)$ is an isotopy invariant.

The skein relation of $H(K)$ can be derived directly, when braiding matrix $\mathcal{B}$ is known:

\begin{figure}[tbp]
\centering
\includegraphics[width=.8\textwidth,trim=0 200 0 200,clip]{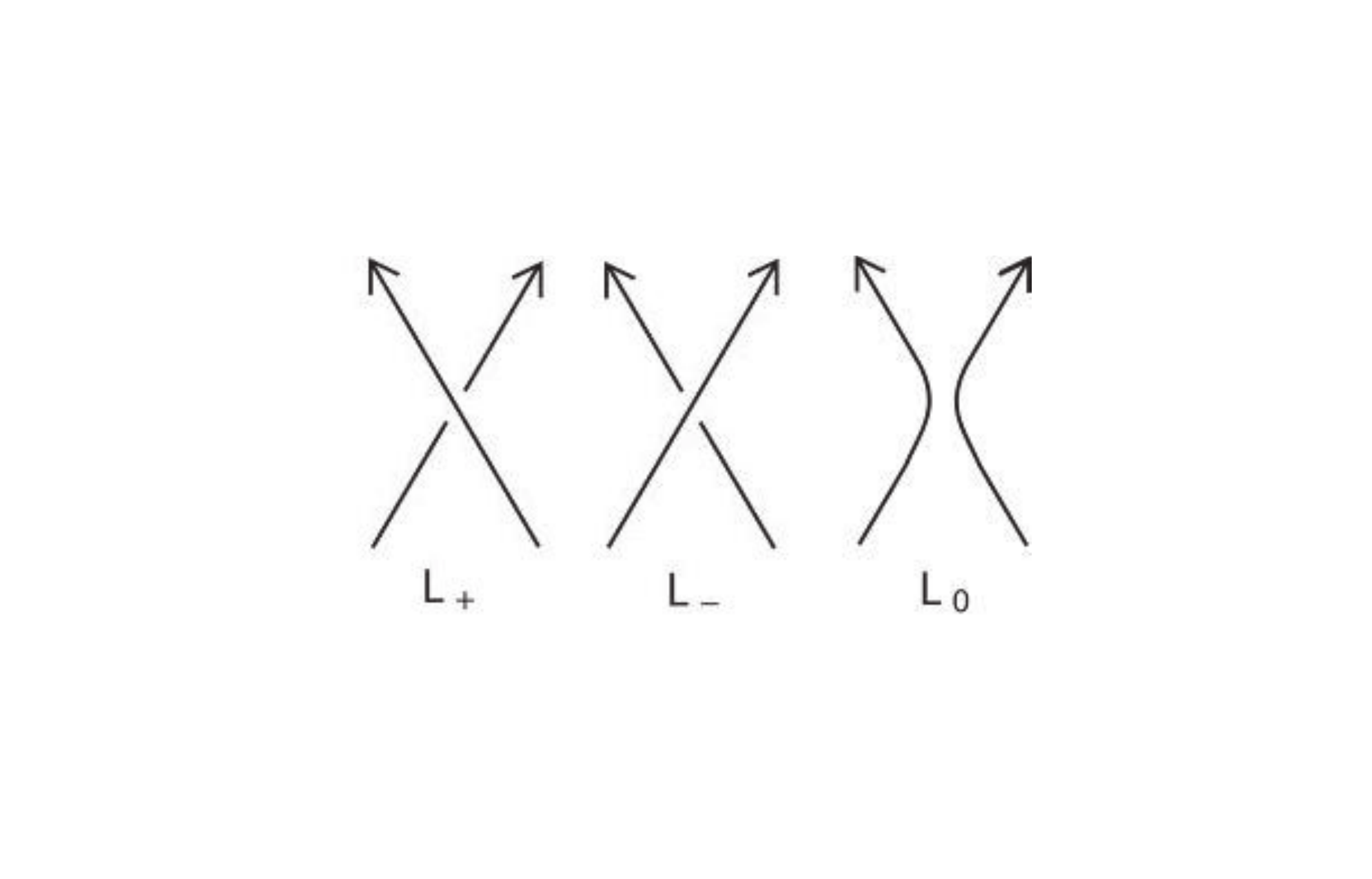}
\caption{\label{fig:skein}}
\end{figure}
\begin{theorem}For the fundamental representation of $A_{n}$ Lie algebra, associated knots invariant $H(K)$ satisfies following skein relation:
\begin{equation}
q^{\frac{n+1}{2}}H(L_{+})-q^{-\frac{n+1}{2}}H(L_{-})=(q^{\frac{1}{2}}-q^{-\frac{1}{2}})H(L_{0}).
\end{equation}
\end{theorem}
Proof: From Theorem \ref{Thm:braiding},
\begin{equation}
q^{-\frac{1}{2(n+1)}}\mathcal{B}-q^{\frac{1}{2(n+1)}}\mathcal{B}^{-1}=(q^{-\frac{1}{2}}-q^{\frac{1}{2}})\mathcal{I}.
\end{equation}
For oriented crossing, $\mathcal{B}$ and $\mathcal{B}^{-1}$ are associated to $<L_{-}>$ and $<L_{+}>$ respectively.
Thus, we have skein relation of $<K>$: \begin{equation}
q^{\frac{1}{2(n+1)}}<L_{+}>-q^{-\frac{1}{2(n+1)}}<L_{-}>=(q^{\frac{1}{2}}-q^{-\frac{1}{2}})<L_{0}>.
\end{equation}

Replacing $<L_{+}>$ and $<L_{-}>$ by $q^{\frac{n(n+2)}{2(n+1)}}H(L_{+})$ and $q^{-\frac{n(n+2)}{2(n+1)}}H(L_{-})$, we have

\begin{equation}
q^{\frac{n+1}{2}}H(L_{+})-q^{-\frac{n+1}{2}}H(L_{-})=(q^{\frac{1}{2}}-q^{-\frac{1}{2}})H(L_{0}).
\end{equation}

This completes the proof.

This skein relation shows that the isotopy invariant $H(K)$ from the generalized Yang-Yang function associated with the fundamental representation of type $A_{n}$ Lie algebra is HOMFLY polynomial ~\cite{c}.

\acknowledgments

 We would like to thank A. Losev, E. Witten, K. Wu and W.-L. Yang for very helpful discussions. Comments and
discussions with Xuexing Lu, Kaiwen Sun, Xiaoyu Jia and Yongjie Wang are also of great help.
This work is partially supported by the National Natural Science Foundation of Grant number
11031005, the Kavli Institute for Theoretical Physics China at the Chinese Academy of Sciences, the School of Mathematical Sciences at Capital Normal University and the Wu
Wen Tsun Key Lab of Mathematics of Chinese Academy of Sciences at University of Science and
Technology of China.

\section{Appendix}

In the following appendix, we prove that, for an unoriented knot $K$, $<K>$ is a regular isotropy invariant, i.e. $<K>$ is invariant under Reidemeister moves of type II and type III.

\subsection*{ $<K>$ is invariant under Reidemeister move of type II }
\begin{figure}[tbp]
\centering
\includegraphics[width=.8\textwidth,trim=0 200 0 200,clip]{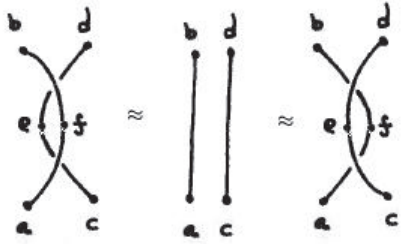}
\caption{\label{fig:VRM} Vertical Reidemeister move II}
\end{figure}

Since the operator $\mathcal{B}$ and $\mathcal{B}^{-1}$ are inverse to each other, i.e.
\begin{equation}
\sum_{e,f}(\mathcal{B}^{-1})^{ef}_{ac}\mathcal{B}^{bd}_{ef}=\delta_{a}^{b}\delta_{c}^{d}=\sum_{e,f}\mathcal{B}^{ef}_{ac}(\mathcal{B}^{-1})^{bd}_{ef},
\end{equation}
$<K>$ is clearly invariant under the vertical Reidemeister move II (see Figure \ref{fig:VRM}).

\begin{figure}[tbp]
\centering
\includegraphics[width=.8\textwidth,trim=0 200 0 200,clip]{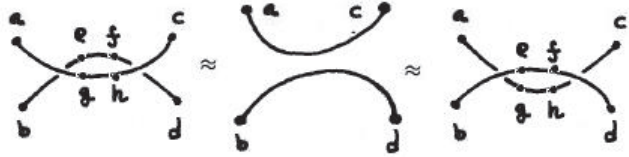}
\caption{\label{fig:HRM} Horizontal Reidemeister move II}
\end{figure}

For horizontal Reidemeister move II (see Figure \ref{fig:HRM}), we have following equation:

\begin{proposition}

\begin{equation}\label{equ:RM2up}
\sum_{e,f,g,h}\mathcal{B}^{ae}_{bg}(\mathcal{B}^{-1})^{fc}_{hd}\mathcal{M}_{ef}\mathcal{M}^{gh}=\mathcal{M}_{bd}\mathcal{M}^{ac}
\end{equation}

\begin{equation}\label{RM2dn}
\sum_{e,f,g,h}(\mathcal{B}^{-1})^{ae}_{bg}\mathcal{B}^{fc}_{hd}\mathcal{M}_{ef}\mathcal{M}^{gh}=\mathcal{M}_{bd}\mathcal{M}^{ac}
\end{equation}

\end{proposition}
Proof: We only prove the equation (\ref{equ:RM2up}), the proof of (\ref{RM2dn}) is similar.
First, consider the left hand side case by case.
\begin{description}
        \item[$c\neq n-a$] \ \begin{itemize}
                             \item if $a\neq b$, $$\hbox{L.H.S}=\mathcal{B}^{ab}_{ba}(\mathcal{B}^{-1})^{n-bc}_{n-ad}\mathcal{M}_{b,n-b}\mathcal{M}^{a,n-a}.$$ From $c\neq n-a$ and $a\neq b$,  $(\mathcal{B}^{-1})^{n-bc}_{n-ad}=0$. Thus, L.H.S=0.
                             \item if $a=b$ and $c\neq d$, $$\hbox{L.H.S}=\mathcal{B}^{aa}_{aa}(\mathcal{B}^{-1})^{n-ac}_{n-ad}\mathcal{M}_{a,n-a}\mathcal{M}^{a,n-a}=0.$$
                             \item if $a=b$, $c=d$, $c>n-a$, then $a>n-c$. $$\hbox{L.H.S}=\sum_{e}\mathcal{B}^{ae}_{ae}(\mathcal{B}^{-1})^{n-ec}_{n-ec}\mathcal{M}_{e,n-e}\mathcal{M}^{e,n-e}=0.$$
                             \item if $a=b$, $c=d$, $c<n-a$, then $a<n-c$.
                                 \begin{equation}
                                 \begin{split}
     &\quad \hbox{L.H.S}\\
     &=\sum_{a\leq e\leq n-c}\mathcal{B}^{ae}_{ae}(\mathcal{B}^{-1})^{n-ec}_{n-ec}\mathcal{M}_{e,n-e}\mathcal{M}^{e,n-e}\\
     &=\mathcal{B}^{aa}_{aa}(\mathcal{B}^{-1})^{n-ac}_{n-ac}\mathcal{M}_{a,n-a}\mathcal{M}^{a,n-a}+\sum_{a< e< n-c}\mathcal{B}^{ae}_{ae}(\mathcal{B}^{-1})^{n-ec}_{n-ec}\mathcal{M}_{e,n-e}\mathcal{M}^{e,n-e}+\mathcal{B}^{an-c}_{an-c}(\mathcal{B}^{-1})^{cc}_{cc}\mathcal{M}_{n-c,c}\mathcal{M}^{n-c,c}\\
     &=(1-q^{-1})\mathcal{M}_{a,n-a}\mathcal{M}^{a,n-a}+(1-q)(1-q^{-1})\sum_{a<e<n-c}\mathcal{M}_{e,n-e}\mathcal{M}^{e,n-e}+(1-q)\mathcal{M}_{n-c,c}\mathcal{M}^{n-c,c}\\
     &=(1-q^{-1})q^{n-c-a}\mathcal{M}_{n-c,c}\mathcal{M}^{n-c,c}+(1-q)(1-q^{-1})q\mathcal{M}_{n-c,c}\mathcal{M}^{n-c,c}\frac{1-q^{n-c-a-1}}{1-q}+(1-q)\mathcal{M}_{n-c,c}\mathcal{M}^{n-c,c}\\
     &=(1-q^{-1})\mathcal{M}_{n-c,c}\mathcal{M}^{n-c,c}(q^{n-c-a}+q(1-q^{n-c-a-1}))+(1-q)\mathcal{M}_{n-c,c}\mathcal{M}^{n-c,c}\\
     &=0.
\end{split}
\end{equation}
                           \end{itemize}
        \item[$d\neq n-b$]\ \begin{itemize}
                            \item if $a\neq b$, then $(\mathcal{B}^{-1})^{n-bc}_{n-ad}=0$. Thus, L.H.S=0.
                            \item if $a=b$ and $c\neq d$, $$\hbox{L.H.S}=\sum_{e}\mathcal{B}^{ae}_{ae}(\mathcal{B}^{-1})^{n-ec}_{n-ed}\mathcal{M}_{e,n-e}\mathcal{M}^{e,n-e}=0.$$
                            \item if $a=b$, $c=d$, $c>n-a$, $$\hbox{L.H.S}=\sum_{e}\mathcal{B}^{ae}_{ae}(\mathcal{B}^{-1})^{n-ec}_{n-ec}\mathcal{M}_{e,n-e}\mathcal{M}^{e,n-e}=0.$$
                            \item if $a=b$, $c=d$, $c<n-a$,
$$\hbox{L.H.S}=\sum_{a\leq e\leq n-c}\mathcal{B}^{ae}_{ae}(\mathcal{B}^{-1})^{n-ec}_{n-ec}\mathcal{M}_{e,n-e}\mathcal{M}^{e,n-e}=0.$$
                          \end{itemize}
        \item[$c=n-a,d=n-b$]\begin{itemize}
                              \item if $a\neq b$, $$\hbox{L.H.S}=\mathcal{B}^{ab}_{ba}(\mathcal{B}^{-1})^{n-bn-a}_{n-an-b}\mathcal{M}_{b,n-b}\mathcal{M}^{a,n-a}=\mathcal{M}_{b,n-b}\mathcal{M}^{a,n-a}.$$
                              \item if $a=b$,
                                  $$\hbox{L.H.S}=\mathcal{B}^{aa}_{aa}(\mathcal{B}^{-1})^{n-an-a}_{n-an-a}\mathcal{M}_{a,n-a}\mathcal{M}^{a,n-a}=\mathcal{M}_{a,n-a}\mathcal{M}^{a,n-a}.$$
                            \end{itemize}
      \end{description}
Clearly, the right hand side $\mathcal{M}_{bd}\mathcal{M}^{ac}$ is :
\begin{description}
  \item[$c\neq n-a$] $$\hbox{R.H.S}=0$$
  \item[$d\neq n-b$] $$\hbox{R.H.S}=0$$
  \item[$c=n-a=d=n-b$] $$\hbox{R.H.S}=\mathcal{M}_{a,n-a}\mathcal{M}^{a,n-a}$$
  \item[$c=n-a\neq d=n-b$] $$\hbox{R.H.S}=\mathcal{M}_{b,n-b}\mathcal{M}^{a,n-a}$$
\end{description}
This completes the proof.

\subsection*{   Braiding matrix $\mathcal{B}$ satisfies Yang-Baxter equation }

Let $V_{1}$, $V_{2}$ and $V_{3}$ be three fundamental representation spaces for $A_{n}$ Lie algebra associated to three parameters $z_{1}$, $z_{2}$ and $z_{3}$. The braiding matrix $\mathcal{B}$ is defined on the tensor product of two representation spaces. We define $$\mathcal{B}_{12}:V_{1}\otimes V_{2}\otimes V_{3}\longrightarrow V_{1}\otimes V_{2}\otimes V_{3}$$ \begin{equation}\mathcal{B}_{12}=\mathcal{B}\otimes id_{3},\end{equation}
and
$$\mathcal{B}_{23}:V_{1}\otimes V_{2}\otimes V_{3}\longrightarrow V_{1}\otimes V_{2}\otimes V_{3}$$
\begin{equation}\mathcal{B}_{23}=id_{1}\otimes\mathcal{B}.\end{equation}
Let $\mathcal{J}_{a,b,c}$ be the vector in the space $V_{1}\otimes V_{2}\otimes V_{3}$, $a,b,c=0,1,2,...,n$. We have checked that
\begin{equation}\mathcal{B}_{12}\mathcal{B}_{23}\mathcal{B}_{12}\mathcal{J}_{a,b,c}=\mathcal{B}_{23}\mathcal{B}_{12}\mathcal{B}_{23}\mathcal{J}_{a,b,c}\hbox{ for all a,b,c}.\end{equation}

From Theorem \ref{Thm:braiding}, we know there are three cases in the braiding. Thus, we just need to calculate 13 cases as follows.
$\mathcal{B}_{12}\mathcal{B}_{23}\mathcal{B}_{12}\mathcal{J}_{a,b,c}$ and $\mathcal{B}_{23}\mathcal{B}_{12}\mathcal{B}_{23}\mathcal{J}_{a,b,c}$ are equal and they are
\begin{description}
  \item[$a>b>c$] $$\gamma^{3}\mathcal{J}_{c,b,a};$$
  \item[$a>b=c$] $$\alpha\gamma^{2}\mathcal{J}_{b,c,a};$$
  \item[$a>c>b$] $$\gamma^{3}\mathcal{J}_{c,b,a}+(\alpha-\beta)\gamma^{2}\mathcal{J}_{b,c,a};$$
  \item[$a=c>b$] $$\alpha\gamma^{2}\mathcal{J}_{a,b,c}+(\alpha-\beta)\alpha\gamma\mathcal{J}_{b,a,c};$$
  \item[$c>a>b$] $$\gamma^{3}\mathcal{J}_{c,b,a}+(\alpha-\beta)\gamma^{2}\mathcal{J}_{b,c,a}+\gamma^{2}(\alpha-\beta)\mathcal{J}_{a,b,c}+(\alpha-\beta)^{2}\gamma\mathcal{J}_{b,a,c};$$
  \item[$a=b>c$] $$\alpha\gamma^{2}\mathcal{J}_{c,a,b};$$
  \item[$a=b=c$] $$\alpha^{3}\mathcal{J}_{a,b,c};$$
  \item[$c>a=b$] $$\alpha\gamma^{2}\mathcal{J}_{c,a,b}+(\alpha-\beta)\alpha\gamma\mathcal{J}_{a,c,b}+\alpha^{2}(\alpha-\beta)\mathcal{J}_{a,b,c};$$
  \item[$a<b<c$] $$\gamma^{3}\mathcal{J}_{c,b,a}+(\alpha-\beta)\gamma^{2}\mathcal{J}_{b,c,a}+\gamma^{2}(\alpha-\beta)\mathcal{J}_{a,b,c}+(\alpha-\beta)\gamma^{2}\mathcal{J}_{c,a,b}+(\alpha-\beta)^{2}\gamma\mathcal{J}_{a,c,b}+\gamma(\alpha-\beta)^{2}\mathcal{J}_{b,a,c}+(\alpha-\beta)^{3}\mathcal{J}_{a,b,c};$$
  \item[$a<b=c$] $$\alpha\gamma^{2}\mathcal{J}_{b,c,a}+(\alpha-\beta)\gamma^{2}\mathcal{J}_{a,b,c}+\alpha\gamma(\alpha-\beta)\mathcal{J}_{b,a,c}+(\alpha-\beta)^{2}\alpha\mathcal{J}_{a,b,c};$$
  \item[$a<c<b$] $$\gamma^{3}\mathcal{J}_{c,b,a}+(\alpha-\beta)\gamma^{2}\mathcal{J}_{a,b,c}+\gamma^{2}(\alpha-\beta)\mathcal{J}_{c,a,b}+(\alpha-\beta)^{2}\gamma\mathcal{J}_{a,c,b};$$
  \item[$a=c<b$] $$\alpha\gamma^{2}\mathcal{J}_{a,b,c}+\alpha\gamma(\alpha-\beta)\mathcal{J}_{a,c,b};$$
  \item[$c<a<b$] $$\gamma^{3}\mathcal{J}_{c,b,a}+\gamma^{2}(\alpha-\beta)\mathcal{J}_{c,a,b},$$
\end{description}where $\alpha=q^{-\frac{n}{2(n+1)}}$, $\beta=q^{\frac{n+2}{2(n+1)}}$ and $\gamma=q^{\frac{1}{2(n+1)}}$ satisfying $\gamma^{2}=\alpha\beta$.

Thus, $\mathcal{B}_{12}$ and $\mathcal{B}_{23}$ satisfy Yang-Baxter equation:
\begin{equation}\mathcal{B}_{12}\mathcal{B}_{23}\mathcal{B}_{12}=\mathcal{B}_{23}\mathcal{B}_{12}\mathcal{B}_{23}.\end{equation}

Similarly, we can prove that $\mathcal{B}^{-1}$ also satisfies Yang-Baxter equation:

\begin{equation}\mathcal{B}^{-1}_{12}\mathcal{B}^{-1}_{23}\mathcal{B}^{-1}_{12}=\mathcal{B}^{-1}_{23}\mathcal{B}^{-1}_{12}\mathcal{B}^{-1}_{23}.\end{equation}

\begin{figure}[tbp]
\centering
\includegraphics[width=.8\textwidth,trim=0 200 0 200,clip]{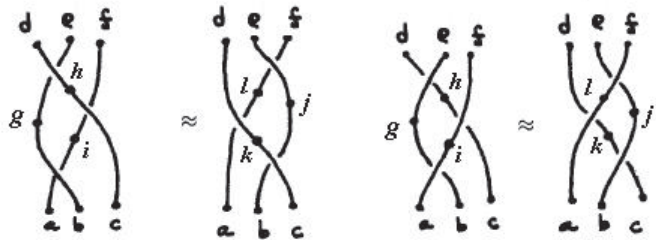}
\caption{\label{fig:YB} Reidemeister move III}
\end{figure}

Write Yang-Baxter equations for $\mathcal{B}$ and $\mathcal{B}^{-1}$ into the contraction of the tensor (see Figure \ref{fig:YB}):
\begin{proposition}

\begin{equation}
\sum_{g,h,i}\mathcal{B}^{gi}_{ab}\mathcal{B}^{hf}_{ic}\mathcal{B}^{de}_{gh}=\sum_{j,k,l}\mathcal{B}^{kj}_{bc}\mathcal{B}^{dl}_{ak}\mathcal{B}^{ef}_{lj}
\end{equation}

\begin{equation}
\sum_{g,h,i}(\mathcal{B}^{-1})^{gi}_{ab}(\mathcal{B}^{-1})^{hf}_{ic}(\mathcal{B}^{-1})^{de}_{gh}=\sum_{j,k,l}(\mathcal{B}^{-1})^{kj}_{bc}(\mathcal{B}^{-1})^{dl}_{ak}(\mathcal{B}^{-1})^{ef}_{lj}
\end{equation}

\end{proposition}

This means that $<K>$ is invariant under Reidemeister move of type III.


\begin{thebibliography}{99}

\bibitem{a}
Jones V F R, \emph{A POLYNOMIAL INVARIANT FOR KNOTS VIA VON NEUMANN ALGEBRAS}, (1985).

\bibitem{b}
Freyd P, Yetter D, Hoste J, et al. \emph{A new polynomial invariant of knots and links},  \emph{Bulletin of the American Mathematical Society} {\bf vol 12.2} (1985) 239-246.


\bibitem{c}
Kauffman, Louis H, \emph{Knots and physics}, {\bf vol 1} (1991).

\bibitem{d}
Witten, Edward, \emph{Quantum field theory and the Jones polynomial}, \emph{Communications in Mathematical Physics} {\bf vol 121.3} (1989) 351-399.

\bibitem{e}
Gaiotto, Davide, and Edward Witten, \emph{Knot invariants from four-dimensional gauge theory}, \emph{Advances in Theoretical and Mathematical Physics} {\bf vol 16.3 } (2012) 935-1086.

\bibitem{f}
 Feigin, Boris, Edward Frenkel, and Nikolai Reshetikhin, \emph{Gaudin model, Bethe ansatz and critical level}, \emph{Communications in Mathematical Physics} {\bf vol 166.1} (1994) 27-62.

\bibitem{g}
 Awata, Hidetoshi, Akihiro Tsuchiya, and Yasuhiko Yamada, \emph{Integral formulas for the WZNW correlation functions}, \emph{Nuclear Physics B} {\bf vol 365.3} (1991) 680-696.

\bibitem{h}
Witten, Edward, \emph{Analytic continuation of Chern-Simons theory}, \emph{Chern-Simons Gauge Theory} {\bf vol 20} (2010) 347-446.

\bibitem{i}
Frenkel, Edward, \emph{Free field realizations in representation theory and conformal field theory}, \emph{Proceedings of the International Congress of Mathematicians} 1995.

\bibitem{j}
Fan, Huijun, \emph{Schrodinger equations,deformation theory and $t t^{*}$ geometry}, arXiv: 1107.1290.

\bibitem{k} Losev, Andrey. \emph{"Hodge strings" and elements of K. Saito¡¯s theory of the Primitive form}, arXiv:hep-th/9801179.

\bibitem{l} Hu, Sen and Liu, Peng. \emph{Knot invariants from the Yang-Yang function,} \emph{to appear in the 6th ICCM Proceedings}.

\bibitem{m} Hu, Sen and Liu, Peng, \emph{Kauffman polynomials from the Yang-Yang function}, \emph{in preparation}.

\end{thebibliography}
\end{document}